\documentclass[reqno,centertags, 12pt]{amsart}

\usepackage{amssymb,amsmath,amsfonts,amssymb} 
-\marginparwidth \oddsidemargin -4mm \evensidemargin
\oddsidemargin
\usepackage{latexsym}
\advance\hoffset by 5mm

\def\@abssec#1{\vspace{.05in}\footnotesize \parindent .2in
{\bf #1. }\ignorespaces}
\newtheorem{theorem}{Theorem}[section]

\newtheorem{lemma}[theorem]{Lemma}
\newtheorem{proposition}[theorem]{Proposition}
\newtheorem{corollary}[theorem]{Corollary}

\def \Rm {\mathbb R}

\def \Zm {\mathbb Z}
\newcommand{\eps}{\varepsilon}

\DeclareMathOperator{\sgn}{sgn}

\newcommand{\lb}{\label}
\newcommand{\la}{\lambda}
\allowdisplaybreaks \numberwithin{equation}{section}

\title{On Discrete Models of the Euler Equation}
\author{Alexander Kiselev and Andrej Zlato\v s}
\thanks{Department of
Mathematics, University of Wisconsin, Madison, WI 53706, USA;
e-mail: kiselev@math.wisc.edu, zlatos@math.wisc.edu}

\begin{document}

\begin{abstract}
We consider two discrete models for the Euler equation describing
incompressible fluid dynamics. These models are infinite coupled
systems of ODEs for the functions $u_j$ which can be thought of as
wavelet coefficients of the fluid velocity. The first model has been
proposed and studied by Katz and Pavlovi\'c. The second has been
recently discussed by Waleffe and goes back to Obukhov studies of
the energy cascade in developed turbulence. These are the only basic
models of this type satisfying some natural scaling and conservation
conditions.
We prove that the Katz-Pavlovi\'c model leads to finite time blowup
for any initial datum, while the Obukhov model has a global solution
for
any sufficiently smooth initial datum.
\end{abstract}

\maketitle

\section{Introduction}\label{intro}

The regularity of solutions to the incompressible Euler equation
in dimension three remains one of the most important open problems
of mathematical fluid dynamics. Recently, a number of simpler
models have been proposed and studied by several authors as a way
to gain insight into the possible behavior of solutions to Euler
and Navier-Stokes equations. Different models have been suggested
by Katz and Pavlovi\'c \cite{KP2}, Friedlander and Pavlovi\'c
\cite{FP}, Dinaburg and Sinai \cite{DS} and Waleffe \cite{W}.
Although these models are fairly drastic simplifications of the
original problem, they do keep a few of the most important
characteristic features of Euler equations. Moreover, we will
argue below that some of these models are quite natural in their
own right as they constitute the simplest class satisfying certain
scaling and dimensional conditions.

A model proposed by Katz and Pavlovi\'c \cite{KP2} is based,
formally, on a wavelet expansion of a scalar function $u(x,t),$ $x
\in \Rm^3,$ over a set of dyadic cubes in $\Rm^3.$ The dyadic cubes
are cubes with the side lengths $2^j,$ $j \in \Zm,$ with vertices at
the points of $2^j \Zm^3.$ If $Q$ is a dyadic cube of size $2^j,$
then its parent $\tilde{Q}$ is a cube with side length $2^{j+1}$
containing $Q.$ Define $C^1(Q)$ the set of all $8$ children of $Q,$
each having side length $2^{j-1},$ and more generally $C^m(Q)$ the
set of all $2^{3m}$ $m^{\rm th}$ generation ``descendants'' of $Q$.
The Katz-Pavlovi\'c model equations describing the evolution of the
wavelet coefficient of $u(x,t)$ corresponding to the cube $Q$ are
the given by \cite{KP2}
\begin{equation}\label{trueKP}
\frac{d u_Q}{dt} = 2^{5j/2} u_{\tilde{Q}}^2 - 2^{5(j+1)/2} u_Q
\sum\limits_{Q' \in C^1(Q)} u_{Q'}.
\end{equation}
The model has quadratic nonlinearity and (formally) conserves the
energy $\sum_{Q} u_Q(t)^2.$ It has been motivated to some extent by
the work \cite{KP1}, where partial regularity of the weak solutions
to the Navier-Stokes equations with hyperdissipation was studied.
The approach of \cite{KP1} is based on controlling the "wavelet
coefficients" of the solution $u_Q=\|\phi_{Q}(x) P_j u\|_2,$ where
$P_j$ are Littlewood-Paley projections restricting the Fourier
transform $\hat{u}(\xi)$ to the annulus of size $\sim 2^j,$ and
$\phi_{Q}$ is a certain smooth function supported on a cube $Q$ of
size $2^{-j(1-\epsilon)},$ $\epsilon
>0.$ The coupled system one gets for the wavelet coefficients from
the Navier-Stokes (or, in our case, Euler) equations is complex, and
\eqref{trueKP} can be obtained from it by dropping all but a few
terms. Thus, $u_Q$ can be roughly thought of as "wavelet
coefficients" describing parts of the solution localized in the cube
$Q$ and in the Fourier space at about $|\xi| \sim 2^j.$ The choice
of the scaling factors in \eqref{trueKP} is determined by the
relation $\|w_Q\|_{\infty} \sim 2^{3j/2} \|w_Q\|_2$ for a wavelet
$w_Q$ supported on a dyadic cube $Q$ of side length $2^j$ in $\Rm^3$
and the bound $\|(u \cdot \nabla)u\|_2 \leq \|u\|_\infty \|\nabla
u\|_2$ (see \cite{KP1,KP2} for more details).

In \cite{KP2} Katz and Pavlovi\'c showed, in particular, that for
any $\eps>0$, there exist initial data $u_{j}(0) \in
H^{3/2+\epsilon}$ which lead to blowup in a finite time. Friedlander
and Pavlovi\'c \cite{FP} considered a related vector model where
they also prove blowup in a finite time. Recently, Waleffe \cite{W}
proposed a simplified model which instead of the branching structure
of the coupled coefficients constitutes a linear tree of the
functions $u_j(t)$ satisfying an infinite system of differential
equations
\begin{equation}\label{KP}
u'_j = \lambda^{j} u_{j-1}^2 - \lambda^{j+1} u_j u_{j+1},
\,\,\,j>j_0, \,\,\,u'_{j_0}= -\lambda^{j_0+1} u_{j_0} u_{j_0+1}.
\end{equation}
Here $\lambda>1$ is a parameter, and $j_0$ is an index corresponding
to the largest relevant space scale (for instance, a period in the
periodic setting). Without loss of generality, we will set $j_0=0$
for the rest of the paper. The original Katz-Pavlovi\'c model
reduces to the system \eqref{KP} with $\lambda=2$ if one assumes
that the coefficients of all cubes of the same side length are the
same.
It is natural to define the Sobolev spaces associated with
\eqref{KP} as
\[
H^s := \{ u_j \, \big| \, \|\{u_j\}\|^2_{H^s}\equiv \sum_{j \geq
j_0} \lambda^{2sj} |u_j|^2 < \infty \}.
\]
Waleffe proved that there exist initial data for which the blowup in
\eqref{KP} happens in any $H^s,$ $s>0,$ and suggested a different
model, given by
\begin{equation}\label{O}
u'_j = \lambda^{j} u_{j-1}u_j - \lambda^{j+1} u^2_{j+1},
\,\,\,j>0, \,\,\,u'_{0}= -\lambda u^2_{1}.
\end{equation}
This model goes back to the work of Obukhov \cite{Ob} who proposed
it in a paper devoted to atmosphere studies
as a simple model for studying the cascade mechanism of energy
transfer in the developed turbulence.
It has been shown in
\cite{W} that the model \eqref{KP} may be related to the inviscid
Burger's equation, making blowup not surprising. In particular,
this model has a built in mechanism of transferring the energy to
higher modes.
On the other hand, the Obukhov model lacks this mechanism and is
thus more subtle and perhaps more realistic. Moreover, in
Proposition~\ref{allmodels} we prove that these models constitute
two basic building blocks of all linear tree coupled mode models
satisfying four natural conditions: a quadratic nonlinearity,
appropriate scaling corresponding to the $(u \cdot \nabla)u$ term,
energy conservation, and nearest neighbor coupling. All of these
except the last one are the features derived from the Euler
equation; the last condition is clearly a simplification designed to
make the problem tractable. Our main goal in this note is to prove
the following two theorems, which to some extent confirm the above
sentiment. For the rest of the paper we call, following Waleffe,
model \eqref{KP} the KP model and model \eqref{O} Obukhov model.

\begin{theorem}\label{thm1}
In the KP model, any non-zero initial datum belonging to $H^1$ leads
to a finite time blowup (in $H^1$).
\end{theorem}

We note that the $H^1$ condition is needed in general to show local
existence of solutions; we discuss this point in
section~\ref{prelim}. If one accepts a parallel between the KP model
and inviscid Burger's equation, the result is not surprising.
Indeed, any non-constant initial datum for the Burger's equation
with periodic boundary conditions leads to blowup in finite time.

On the other hand, solutions of the Obukhov model are regular.

\begin{theorem}\label{thm2}
In the Obukhov model, the solution corresponding to any initial
datum in $H^s,$ $s>1,$ is regular for all times. That is, for any
$u_0 \in H^s$ with $s>1$ and for any $T>0$ there exists a unique
solution $\{ u_j \} \in C([0,T], H^s)$ such that $u_j(0)=(u_0)_j.$
\end{theorem}


This theorem is probably the most interesting, and certainly the
most subtle and difficult to prove result of this paper. It
demonstrates an intriguing dichotomy between the properties of two
basic dyadic models.

For generic initial data in the Obukhov model, we have a stronger
regularity and even dissipation properties, in the following sense.

\begin{theorem}\label{thm3}
Let $b_j(\omega)$ be independent uniformly bounded random variables
such that the probability of $b_j(\omega)$ being nonpositive is
uniformly bounded away from zero: $P[b_j(\omega)\le 0]>\rho>0.$
Assume $a_j>0$ are such that $\sum_j \lambda^{2sj}|a_j|^2 < \infty,$
$s > 1.$ Then with probability one a solution $\{u_j(t)\}$ of the
Obukhov model corresponding to the initial datum $u_j(0)=a_j
b_j(\omega)$ satisfies $\|u\|_{H^r} \leq C(r,\omega)$ for all times
$t$ and any $r<s.$ Moreover, as $t \rightarrow \infty,$ we have
\begin{equation}\label{sobas}
\lim_{t\to\infty} \|u(t)\|^2_{H^r} = \lim_{t\to\infty} u_0(t)^2 =
E_0 \equiv \sum\limits_{j\ge 0} u_j(0)^2,
\end{equation}
that is, the solution $u$ converges in $H^r$ to a constant solution
with all energy concentrated in the lowest mode.
\end{theorem}

We describe some finer properties of the dynamics of the KP and
Obukhov models as well.

There are many interesting questions that remain open. In
particular, whether Theorem~\ref{thm2} holds for $s=1$. Other
natural questions include global existence of solutions in the
branched Obukhov model (an analog of \eqref{trueKP}) and in the
Navier-Stokes version of \eqref{O}. It seems reasonable to expect
that regularity results for \eqref{O} should carry over to these
cases. Clearly, the analog of the Laplacian term only adds
dissipation, and branching is likely to make energy cascade towards
high level modes harder to realize. However, on the technical level,
the questions are not trivial due to the subtleties of the proof of
Theorem~\ref{thm2}. We did not attempt to address these issues here
to keep the present paper from becoming overly technical.

We note that models similar in spirit to \eqref{KP} and \eqref{O}
--- shell models --- have been studied in the physics literature for a
long time (see, e.g., \cite{Gle,Fri,Okh,Sabra}, and \cite{Bif} for a
recent review). One version of these models, the "Sabra" shell
model, has recently been studied analytically in \cite{CLT}. What
makes questions like existence and (in some sense) regularity of
solutions easier to treat in the shell models setting than in
\eqref{KP}, \eqref{O} is a weaker scaling factor in the equations
(corresponding, generally speaking, to the scaling assumption $\|(u
\cdot \nabla) u\|_2 \sim \|u\|_2\|\nabla u\|_2$). This leads to the
shell models being "subcritical", that is the nonlinearity is
controlled by the dissipation term. However the models \eqref{KP},
\eqref{O}, even when a term representing Laplacian with appropriate
scaling is added, are "supercritical". It is only certain
monotonicity properties of these models and detailed analysis of
their dynamics that make answering the basic regularity/blowup
questions possible. Many of the subtler results established for the
"Sabra" model in \cite{CLT} appear harder to establish for the
dyadic models of Navier-Stokes equations at this time.

In the next section we collect some preliminary results, postponing
the proof of local existence in $H^1$ of solutions to our models to
an appendix. The proofs of our main theorems appear in Sections
3--5.

\section{Preliminaries}\label{prelim}

In this section we collect and prove some simple useful facts
about the KP and Obukhov models. Let us start by stating the
result on local existence of solutions.

\begin{proposition}\label{P.exist}
Assume that the initial datum $u_j(0)$ for either KP or Obukhov
model lies in $H^s$ for some $s \geq 1.$ Then there exists  a unique
solution $u \in C([0,T], H^s),$ for some time
$T=T(\|u(0)\|_{H^s})>0.$ The $H^s$ norm of this solution satisfies
\begin{equation}\label{exist}
\|u(t)\|_{H^s} \leq \|u(0)\|_{H^s} e^{C \int_0^t {\rm sup}_j
\{\lambda^j u_j(r)\}\, dr}.
\end{equation}
In particular, the solution blows up in finite time $\tau$ only if
$\int_0^\tau {\rm sup}_j \{\lambda^j u_j(r)\}\, dr = \infty.$
\end{proposition}

\begin{proof}
Local existence of solutions has been proved in \cite{FP} using
fixed point arguments. The argument in \cite{FP} is given for the
case of KP model with a specific choice of $\lambda$ (hence our
$H^1$ notation corresponds to $H^{5/2}$ in their setting), but it
can be adapted easily to the Obukhov model as well. We sketch this
argument in the Appendix. Therefore, here we will only discuss
\eqref{exist}. Carrying out the differentiation and substituting the
expression for the time derivatives from \eqref{KP} (resp.
\eqref{O}) we find
\[ \frac{d}{dt} \sum_j \lambda^{2sj} u_j^2(t) \leq
C {\rm sup}_j \{\lambda^j u_j(t)\} \,\sum_{j=0}^\infty \lambda^{2sj}
u_j^2(t),
\] providing the required bound.
\end{proof}

Now we make a few critical observations on the monotonicity
properties of our models. From now on, all properties are stated for
the solutions described in Proposition~\ref{P.exist}, and hold on
the existence interval described in that proposition.

\begin{proposition}\label{keyprop}
The following properties hold for KP and Obukhov models.
\begin{itemize} \item Both KP and Obukhov models
conserve the energy $E_0\equiv\sum_{j \geq 0} |u_j(t)|^2.$ \item In
the KP model, if $u_j(t_0) \geq 0$ for some $t_0,$ then $u_j(t) \geq
0$ for all times $t \geq t_0.$ \item In the Obukhov model, if
$u_j(t_0) \leq 0$  for some $t_0,$ then $u_j(t) \leq 0$ for all
times $t \geq t_0.$
\end{itemize}
\end{proposition}
\begin{proof}

The first property is checked directly by differentiating the
energy. Clearly each $u_j(t)$ is differentiable, and the fact that
solution is $H^1$ allows us to sum the right hand side, obtaining
zero. To prove the last two properties, one just writes explicitly
the expression for $u_j(t).$ For example, in the Obukhov model we
have
\[ u_j(t) = e^{\lambda^j \int_{t_0}^t u_{j-1}(r)\,dr} \left(
u_{j}(t_0) - \lambda^{j+1} \int_{t_0}^t e^{-\lambda^j
\int_{t_0}^\rho u_{j-1}(r)\,dr} u_{j+1}^2(\rho)\,d\rho \right). \]
\end{proof}

Let us define $E_j(t) \equiv \sum_{l \geq j} |u_l(t)|^2.$ Note that
$E_j'(t) = 2\lambda^{j}u_{j-1}^2u_j$ in the KP and $E_j'(t) =
2\lambda^{j}u_{j-1}u_{j}^2$ in the Obukhov model. Hence, in both
models positive coefficients generate energy transfer to higher
modes and negative coefficients transfer energy to lower modes.
Since Proposition~\ref{keyprop} shows that positive coefficients are
stable in the KP model and negative ones are stable in the Obukhov
model, it is not surprising that the latter is more regular. One
more indication of this regularity is the following description of
the dynamics corresponding to initial data with only finite number
of excited modes.

\begin{proposition}\label{fmodes}
In the Obukhov model, if $u_j(0) = 0$ for any $j > j_1,$ then
$u_j(t) =0$ for any $t$ and $j>j_1$. In this case, as time goes to
infinity, all energy concentrates in the first mode $u_{0}.$
Moreover, if $u$ is any solution that remains in $H^1$ for all time,
then $u_{j}(t) \rightarrow 0$ as $t \rightarrow \infty$ for all
$j>0.$
\end{proposition}

\begin{proof}
The first statement is obvious. Let us prove the third statement
(which in turn proves the second in the case of eventually vanishing
$u_j(0)$). It is clear from \eqref{O} that $u_{1}(t) \rightarrow 0,$
or else $u_{0}$ grows unboundedly large negative, contradicting the
energy conservation. This holds since $|u'_{1}(t)| \leq \lambda^2
E_{0},$ so the function $u_1(t)$ cannot just have increasingly
narrow spikes. Now, if $u_j(t) \rightarrow 0,$ then $u_{j+1}(t)
\rightarrow 0$. Otherwise the equation
\[ u'_j(t) = \lambda^j u_{j-1} u_j - \lambda^{j+1} u_{j+1}^2 \]
and $|u_{j-1}|\le\sqrt{E_0}$ give us a contradiction as in the case
$j=0$ above.
\end{proof}

Finally, before proving our main results, we state the following
observation, which is elementary to verify. It shows that the KP
and Obukhov models are basic building blocks of all mode couplings
with certain natural properties.

\begin{proposition}\label{allmodels}
Assume that real valued functions $u_j(t)$ satisfy an infinite
system of differential equations such that:
\begin{itemize}
\item The right hand side is quadratic in $u$ \item The coupling
is nearest neighbor only, that is only $u_{j-1},$ $u_j$ or $u_{j+1}$
may appear in the equation for $u_j'$ \item Each term on the right
hand side of the equation for $u_j$ has a factor of $\lambda^{j}$
times a constant independent of $j$ \item The energy $\sum_j u_j^2$
is conserved.
\end{itemize}
Then the system must have the form
\begin{equation}\label{alsyst}
u_j' = \alpha( \lambda^{j}u_{j-1}^2 - \lambda^{j+1} u_j u_{j+1}) +
\beta (\lambda^j u_j u_{j-1} - \lambda^{j+1} u_{j+1}^2),
\end{equation}
that is, the right hand side must be a linear combination of the KP
and Obukhov models.
\end{proposition}

Theorems \ref{thm1} and \ref{thm2} show that if the initial datum is
in $H^s$, $s>1$, then the solution of \eqref{alsyst} always blows up
resp.~stays regular if $\alpha=1$, $\beta=0$ resp.~$\alpha=1$,
$\beta=0$. It is an interesting open question how the competition of
these two phenomena affects the behavior of solutions of
\eqref{alsyst} when both $\alpha, \beta\neq 0$. Notice that when
$\sgn(\alpha)=\sgn(\beta)$, then we do not have at our disposal a
version of the maximum principle, as are the second and third claims
of Proposition \ref{keyprop}. This structural difference in the
general case will present an extra difficulty in the analysis of the
dynamics of the problem.

\section{Blowup in the Katz-Pavlovi\'c model}

In this section we prove Theorem~\ref{thm1}. We therefore assume,
towards contradiction, that the solution exists in $H^1$ for all
times and $\|u\|_{H^1}$ is locally bounded. Let us define the
``positive'' and ``negative'' energies by
\[ E_{\pm,j}(t) = \sum\limits_{l \geq j, \pm u_l \geq 0} u_l(t)^2. \]
The following lemma shows that for any non-zero initial datum and
any $j$, $E_{+,j}(t)>0$ for $t> t_j.$
\begin{lemma}\label{posturn}
For any non-zero initial datum and any $j>0,$ we have $u_j(t)>0$ for
$t>t_j.$
\end{lemma}
\begin{proof}
Recall that $j_0=0$. Note that
\[
 u_{0}(t) = u_{0}(0) e^{-\lambda \int_0^t
u_{1}(s)\,ds}.
\]
Assume that $u_{1}(0)<0,$ and never turns positive. Then at least
we must have $u_{1}(t) \rightarrow 0$ as $t \rightarrow \infty,$
or else $u_{0}$ grows unbounded. But then we get a contradiction
with the equation
\[ u'_{1}= \lambda u_{0}^2 - \lambda^2 u_{1}
u_{2}, \] since $|u_{0}(t)| \geq |u_{0}(0)|>0$ for all times (if
$u_{0}(0)=0,$ it is never in the play, and so we should start from
$j=1$). Thus $u_{1}$ must become positive. Now if $u_j(t_j)>0,$ then
$u_{j+1}(t)$ must turn positive at some finite time too, by an
argument identical to the above.
\end{proof}

Next, we show that the positive energy is always increasing.

\begin{lemma}\label{posin}
For any $j,$ $E_{+,j}(t)$ is monotone increasing. The negative
energy $E_{-,j}(t)$ is monotone decreasing.
\end{lemma}

\begin{proof}
At any given moment, $E_{+,j}(t)$ can be written as a sum of sums
$\sum_{j_1 \leq l \leq j_2} u_l^2,$ where $u_l(t) \geq 0$ for $j_1
\leq l \leq j_2,$ and $u_{j_1-1}(t),u_{j_2+1}(t) <0$ (or $j_1=0$).
Then
\[
 \frac{d}{dt} \sum_{j_1 \leq l \leq j_2} u_l^2 = 2 \sum_{j_1
\leq l \leq j_2} u_l (\lambda^l u_{l-1}^2 - \lambda^{l+1} u_l
u_{l+1}) = 2(\lambda^{j_1} u_{j_1} u_{j_1-1}^2 - \lambda^{j_2+1}
u_{j_2}^2 u_{j_2+1}) \geq 0.
\]
Moreover, we see from the above
argument that
\begin{equation}\label{poseninc}
E_{+,j}'(t) \geq 2 \lambda^j u_j u_{j-1}^2.
\end{equation}
This bound is not relevant if $u_j(t) \leq 0,$ but we will need it
later in the case when we know that $u_j$ is positive. The proof for
$E_{-,j}$ is similar.
\end{proof}

Theorem~\ref{thm2} will be a simple consequence of the following
key lemma.

\begin{lemma}\label{blowup}
Let $q\in(\lambda^{-1},1)$ and $\rho\equiv (\lambda
q)^{-1}\in(0,1)$, and assume that $j$ is large enough (depending on
$\lambda$, $q$, and $E_0$). Then for any $C>0$ there is
$A=A(C,\la,q)<\infty$ (independent of $j$) so that if $E_{+,j}(t_0)
\geq Cq^j$ for some $t_0$, then there exists a time $t \in
[t_0,t_0+2\tau_j]$, with $\tau_j \equiv A\rho^j$, such that either
$E_{+,j+1}(t) \geq Cq^{j+1}$ or $E_{+,j}(t) \geq 2 Cq^j.$
\end{lemma}


\begin{proof}
Assume that for all $t \in [t_0,t_0+2\tau_j]$ we have $E_{+,j+1}(t)
\leq Cq^{j+1}.$ Then by $E_{+,j}(t)\ge E_{+,j}(t_0)\ge Cq^j,$ we
must have $u_j(t) \geq 0$ and $u^2_j(t) \geq Cq^j(1-q)$ for any $t
\in [t_0,t_0+2\tau_j].$ Let
\[
A\equiv \max \bigg\{ 1, \frac{1+\lambda q}{\sqrt C (1-q)^2},
\frac{4\sqrt E_0}{C(1-q)^2} \bigg\}
\]

Consider first the case where $u_{j+1}(t_1) \geq 0$ for some $t_1
\in [t_0,t_0+\tau_j].$
The amount of energy transfer from $j^{\rm th}$ to $(j+1)^{\rm st}$
mode is bounded from below by (recall \eqref{poseninc})
\[ \int_{t_1}^{t_1+\tau_j} E_{+,j+1}'(t)\,dt
\geq 2\lambda^{j+1} \int_{t_1}^{t_1+\tau_j} u_j(t)^ 2 u_{j+1}(t)
\,dt. \] It must not exceed $Cq^{j+1}$ to avoid contradiction, so
\[ \int_{t_1}^{t_1+\tau_j} u_{j+1}(t)\,dt \leq
\frac{q}{(1-q)\lambda^{j+1}}. \] But
\[ u_{j+1}'(t) = \lambda^{j+1}u_j^2 - \lambda^{j+2} u_{j+1}
u_{j+2}; \] thus
\begin{equation}\label{entranb}
u_{j+1}(t_1+\tau_j)- u_{j+1}(t_1) \geq \lambda^{j+1}Cq^j(1-q)\tau_j
- \lambda^{j+2} \frac{q}{(1-q)\lambda^{j+1}}\sqrt{C} q^{(j+1)/2}.
\end{equation}
The above bound follows from the fact that $u_{j+1}\ge 0$ and
$u_{j+2}\le \sqrt C q^{(j+1)/2}$ on $[t_1,t_1+\tau_j]$, the latter
by our assumption on $E_{+,j+1}$. The right hand side of
\eqref{entranb} equals
\begin{equation}\label{bbc}
 \sqrt{C} q^{(j+1)/2} \left( \sqrt{C} \lambda^{j+1} \tau_j (1-q)
q^{(j-1)/2} - \frac{\lambda q}{1-q}\right). \end{equation}
Since $A\ge (1+\lambda q)/\sqrt C (1-q)^2$ and $\tau_j=A\rho^j$, the
expression in the brackets in \eqref{bbc} is greater than one.

It remains to consider the case where $u_{j+1} (t) < 0$ for $t \in
[t_0,t_0+\tau_j].$ Recall that we have $u_j^2 \geq C q^j(1-q),$ and
$-u_{j+1} u_{j+2} \geq u_{j+1} F_0,$ where $F_0^2=E_0$ is the total
(conserved) energy of the solution. Then from \eqref{KP} we obtain
for any $t_1\in [t_0,t_0+\tau_j]$,
\begin{eqnarray} u_{j+1}(t) \geq u_{j+1}(t_1) e^{\lambda^{j+2}
F_0 (t-t_1)} + C \int_{t_1}^t \lambda^{j+1} q^j (1-q)
e^{\lambda^{j+2} F_0(t-s)}\,ds  \nonumber \\
\geq e^{\lambda^{j+2} F_0(t-t_1)} \left( u_{j+1}(t_1) + C
\lambda^{-1} q^j (1-q)F_0^{-1}(1-e^{-\lambda^{j+2}
F_0(t-t_1)})\right). \label{negb}
\end{eqnarray}
Assume without loss of generality that $j$ is large enough, so that
$ \lambda^j F_0 \rho^j >> 1$ (then also $\lambda^{j+2} F_0\tau_j >>
1$ because $A\ge 1$). If for some $t_1 \in [t_0, t_0+\tau_j/2]$ the
value of $u_{j+1}(t)$ goes above $-\tfrac 12
C\lambda^{-1}q^j(1-q)F_0^{-1},$ we see from \eqref{negb} that
$u_{j+1}(t)$ will become positive before $t_0+\tau_j.$ Thus, we must
have
\[ u_{j+1}(t) \leq -\tfrac 12 C\lambda^{-1}q^j(1-q)F_0^{-1} \]
for $t\in [t_0, t_0+\tau_j/2].$ But then for these $t$,
\[ \frac{d}{dt} u_j^2 \geq -2\lambda^{j+1} u_j^2 u_{j+1} \geq
\lambda^{j}C^2q^{2j}(1-q)^2F_0^{-1}. \] This implies
\[ u_j(t_0+\tau_j/2)^2\ge u_j(t_0+\tau_j/2)^2 - u_j(t_0)^2 \geq \tfrac12
\tau_j\lambda^{j}C^2q^{2j}(1-q)^2F_0^{-1} \geq 2Cq^j \] since $A\ge
4F_0/C(1-q)^2$. Thus, $E_{+,j}(t_0 +\tau_j/2) \geq 2Cq^j,$ and the
lemma is proved.
\end{proof}
The second alternative in Lemma~\ref{blowup} is needed since if
$u_{j+1}$ is very large negative, it seems reasonable that it may
take some time before it becomes positive and the positive energy
starts being transferred up. The proof is based on the observation
that in this case, the negative energy from the $(j+1)^{\rm st}$
mode is quickly transferred into the positive one at the $j^{\rm
th}$ mode. The following corollary shows that actually the lemma
holds in a simpler form, without the second alternative, if we
increase the waiting time slightly.

\begin{corollary}\label{simpler}
In the setting of Lemma \ref{blowup}, there exists $t \in [t_0,
t_0+2 \log_2 (E_0/Cq^j) \tau_j],$ such that $E_{+,j+1}(t) \geq
Cq^{j+1}.$
\end{corollary}
\begin{proof}
Recall that the total energy of the solution is equal to $E_0.$
Applying Lemma~\ref{blowup} repeatedly on the $j^{\rm th}$ level, we
see that the second alternative cannot hold more than $\log_2
(E_0/Cq^j)$ times.
\end{proof}

Now we can complete the proof of Theorem~\ref{thm1}.

\begin{proof}[Proof of Theorem~\ref{thm1}]
Pick some $q\in(\lambda^{-1},1)$ and denote $\tilde{\tau}_j = 2 \log
(E_0/Cq^j) \tau_j.$ It is clear that
\[ \tilde{\tau} = \sum\limits_{j} \tilde{\tau}_j < \infty. \]
Lemma \ref{posturn} shows that each $u_j$ (in particular, those to
which Lemma \ref{blowup} applies) will eventually become positive.
Using Corollary \ref{simpler} one then shows by induction that for
some $t_0<\infty$, $C>0$, and for all large $j$, there exists $t_j
\in [t_0, t_0+\tilde{\tau}]$ such that $E_{+,j}(t) \geq Cq^j$. Note
that $t_j$ can be chosen to be increasing. But then the $H^1$ norm
satisfies
\[ \|u(t_j)\|_{H^1}^2 \geq C\lambda^{2j} q^j \rightarrow \infty
\]
because $q>\lambda^{-1}$. The proof is finished.
\end{proof}

\section{Almost sure estimates in the Obukhov model}

In this section we prove Theorem~\ref{thm3} as a warmup. This result
is rather straightforward, relying only on the fact that negative
coefficients are stable in the Obukhov model and that the energy
always flows to the lower modes across any negative site.

\begin{proof}[Proof of Theorem~\ref{thm3}]
Consider a realization of $\{b_j(\omega)\},$ that has infinitely
many sites $j_1(\omega) < \dots < j_n(\omega) < \dots$ at which
$b_{j_l}(\omega) \leq 0.$ Such realizations occur with probability
1, by the hypothesis. We also set $j_0(\omega)=0$ by convention.
Since
\[ E'_{j_l(\omega)+1}(t) = 2 \lambda^{j_l(\omega)+1}
u_{j_l(\omega)}(t) u_{j_l(\omega)+1}(t)^2, \] we see that by
Proposition~\ref{keyprop}, $E'_{j_l(\omega)+1}(t) \leq 0,$ for all
$t>0$ for which the solution exists. Therefore, for all such times
we have the following estimate
\begin{equation}\label{hrest}
\|u(t)\|_{H^r}^2 \leq \sum\limits_{l=0}^\infty
\lambda^{2j_l(\omega)r} \left( \sum\limits_{m=
j_{l-1}(\omega)}^{j_l(\omega)} |u_m(0)|^2 \right)
\leq C_1 \sum\limits_{l=0}^\infty \lambda^{2j_l(\omega)r
-2j_{l-1}(\omega)s},
\end{equation}
since $|u_m(0)| \leq C \lambda^{-ms}$ by assumption. We claim that
for any $\alpha>0,$ with probability one we have
\begin{equation}\label{bcest}
j_l(\omega) - j_{l-1}(\omega) \leq \alpha j_{l-1}(\omega)
\end{equation}
for all but finitely many $l.$ If that were the case, take $\alpha =
(s-r)/2r.$ Then
\[ 2j_l(\omega)r - 2j_{l-1}(\omega)s \leq -(s-r)j_{l-1}(\omega)
\]
almost surely for all but finitely many $l.$ In that case, the sum
\eqref{hrest} converges almost surely, proving $\|u\|_{H^r}\le
C(r,\omega)$. This and local existence in $H^1$ (note that $1<s$)
now gives the existence of the solution in $H^r$, $r<s$, for all
times.

To prove \eqref{bcest}, split natural numbers into non-overlapping
intervals $L_n \equiv \{ j\,|\,3^{n-1} < j \leq 3^n\}.$ It is clear
that for all $\alpha$ small enough, any interval $I_l=(j_{l-1},j_l)
$ satisfying $j_l-j_{l-1} > \alpha j_{l-1}$ will have an
intersection of size at least $\alpha 3^{n-2}$ with some $L_n$. The
probability of having such an interval of negative $b_j(\omega)$'s
in $L_n$ is less than $3^n (1-\rho)^{\alpha 3^{n-2}}.$ Since the
events of having such an interval in $L_n$ for different $n$ are
independent, we find that the probability of having an infinite
number of such intervals is zero by the Borel-Cantelli lemma.

The fact that $\|u(t)\|_{H^r}^2$ converges to $E_0$ follows from the
above argument and Proposition~\ref{fmodes}. Indeed, with
probability one $u_0(t)^2\le \|u(t)\|_{H^r}^2 \leq
u_0(t)^2+\sum_{l\ge 1}\limits A_l(t),$ where
\[ A_l(t) = \sum\limits_{m= j_{l-1}(\omega)+1}^{j_l(\omega)}
|u_m(t)|^2 \lambda^{2mr}, \] and we saw that $A_l(t) \leq C(\omega)
\lambda^{-(s-r)l}.$ But Proposition~\ref{fmodes} also implies
$A_l(t) \rightarrow 0$ as $t \rightarrow \infty$ for any $l.$ Thus,
by the dominated convergence theorem,
\[
\lim_{t\to\infty} (\|u(t)\|_{H^r}^2 - u_0(t)^2) = 0.
\]
In particular, for $r=0$ we get using energy conservation
\[
E_0 = \|u(t)\|_{L^2}^2 = \lim_{t\to\infty} u_0(t)^2
\]
which yields \eqref{sobas}.
\end{proof}

\section{Regularity in the Obukhov model}

We will now prove Theorem \ref{thm2}. Assume, towards contradiction,
that for some initial datum $u(0)$ with $\|u(0)\|_{H^s} \le 1$ (this
can be assumed without loss of generality, by scaling in $u$ and
$t$), $u$ blows up at time $T<\infty$, that is,
\begin{equation} \lb{5.1}
\limsup_{t\to T} \| u(t) \|_{H^s}=\infty
\end{equation}
and $\| u(t) \|_{H^s}$ is bounded for $t\in[0,T-\eps]$ and any
$\eps>0$ (using \eqref{exist} and $\| u(t) \|_{H^s}\ge \sup_j\{\la^j
u_j(t) \}$, one can actually show that the $\limsup$ must be
$\lim$). Proposition \ref{P.exist} shows that this is only possible
if
\begin{equation} \lb{5.2}
\limsup_{t\to T} \sup_j\{\la^j u_j(t) \}=\infty.
\end{equation}
Although a priori it only follows from the proposition that the
$\limsup$ is $\infty$ for some $T^*\le T$, it is immediate from
$s>1$ that in that case \eqref{5.1} would hold for $T^*$ and so
$T^*=T$. We have
\begin{equation} \lb{5.3}
u_j(0)\le \la^{-sj}
\end{equation}
and by $\|u(t)\|_{L^2}=\|u(0)\|_{L^2}\le \| u(0) \|_{H^s}\le 1$,
\begin{equation} \lb{5.4}
|u_j(t)|\le 1.
\end{equation}
Finally, we recall that
\begin{equation*} 
E_j(t)\equiv \sum_{l\ge j} u_l(t)^2
\end{equation*}
satisfies
\begin{equation} \lb{5.6}
E_j'(t) = 2\la^j u_{j-1}(t)u_j(t)^2
\end{equation}
(with $u_{-1}\equiv 0$).

Our strategy will be to first narrow down the possibility of blowup
to a specific scenario (Lemma \ref{L.5.1}) and then exclude blowup
under this scenario (Lemma \ref{L.5.3}). Let $t_j<T$ be the first
time such that
\begin{equation} \lb{5.7}
u_j(t_j) = \la^{-j}
\end{equation}
(if there is no such time we let $t_j\equiv \infty$). If
$t_j<\infty$, then
\begin{equation} \lb{5.8}
u_j(t)>0 \text{ for } t\in[0,t_j],
\end{equation}
by Proposition \ref{keyprop}. Therefore we can use
\begin{equation} \lb{5.9}
\frac{u_j'}{u_j} = \la^{j} u_{j-1} - \la^{j+1} \frac{u_{j+1}^2}{u_j}
\end{equation}
for $j>0$ to obtain from \eqref{5.3} and \eqref{5.7}
\begin{equation*} 
(s-1)j\log \la \le \log \frac{u_j(t_j)}{u_j(0)} \le \la^j
\int_0^{t_j} u_{j-1}(t)\, dt \le \la^j T \sup_{t\le t_j} u_{j-1}(t).
\end{equation*}
Hence
\begin{equation} \lb{5.11}
\sup_{t\le t_j} u_{j-1}(t)\ge  \frac{(s-1)\log\la}{T\la} j
\la^{-(j-1)},
\end{equation}
which means that $t_{j-1}\le t_j$ once $j>T\la((s-1)\log\la)^{-1}$
(this is obviously true also when $t_j=\infty$). Therefore $t_j$ is
eventually non-decreasing and has a limit $\tau$. Now \eqref{5.2}
and \eqref{5.4} imply that $\tau\neq\infty$ and so $\tau\le T$. Then
\eqref{5.11} shows $t_{j-1}<t_j$ for large $j$ as well as
$\sup_{t\le t_j} u_{j-1}(t)\la^{j-1}\to\infty$ as $j\to\infty$, and
so $T\le\tau$ (since blowup cannot happen before $T$). Hence $t_j$
is eventually increasing and $t_j\to T$. From now on we will
consider $j$ large enough so that $T-1<t_{j}<t_{j+1}<T$ and set
$I_j\equiv [t_{j},t_{j+1}]$. Note that
\begin{equation} \lb{5.11a}
u_j(t_{j+1})\ge 0
\end{equation}
because $t_{j+1}$ is the first time when $u_{j+1}$ reaches
$\la^{-j-1}$, and $u_{j+1}$ (if it is positive) has to decrease when
$u_j<0$. At various places in the argument below we will further
increase the size of $j$ under consideration.

We choose $\eps\in(0,\tfrac{s-1}5)$. For $t\le t_j$ and $l\ge 1$ we
have by \eqref{5.8} and \eqref{5.9},
\[
\log \frac{u_{j+l}(t)}{u_{j+l}(0)} \le \la^{j+l} \int_0^{t}
u_{j+l-1}(\tau)\, d\tau \le \la T
\]
since $t < t_{j+l}$. Therefore by \eqref{5.3},
\begin{equation} \lb{5.12}
u_{j+l}(t) \le e^{\la T} \la^{-s(j+l)} \le \la^{-(s-\eps)(j+l)}
\end{equation}
for large enough $j$, $t\le t_j$, and $l\ge 1$. This and \eqref{5.9}
gives
\begin{equation} \lb{5.13}
(s-1-\eps)j\log\la \le \log \frac{u_{j+1}(t_{j+1})}{u_{j+1}(t_j)}
\le \la^{j+1} \int_{I_j} u_j(t)\, dt.
\end{equation}

Thus for all large $j$, $u_j$ has to become large compared to
$\la^{-j}$ somewhere on $I_j$, while $u_{j+1}$ increases to
$\la^{-j-1}$ and all the higher modes are tiny. This shows that for
blowup at $T$ to occur there must be a ``wave'' of large $\la^ju_j$
moving from low to high modes, reaching infinity in finite time.
Next we will show that this wave has to be eventually very thin.
Namely, we will show that modes just behind the head of the wave
quickly become negative when $j$ is large.

\begin{lemma} \lb{L.5.1}
For all large enough $j$ we have $u_{j-1}(p_j)\le 0$ with $p_j\in
I_j$ defined by
\begin{equation} \lb{5.14}
\la^{j+1} \int_{t_j}^{p_{j}} u_j(t)\, dt = \frac{3(s-1-\eps)}4
j\log\la.
\end{equation}
\end{lemma}

{\it Remark.} Note that this $p_j$ is unique by \eqref{5.13} and
Proposition \ref{keyprop}
\smallskip

Lemma \ref{L.5.1} will be a consequence of the following weaker
formulation of the thin wave property.

\begin{lemma} \lb{L.5.2}
For any $j_1$ there is $j>j_1$ such that $u_{j-2}(r_j)\le 0$ for
 $r_j\in I_j$ defined by
\begin{equation} \lb{5.15}
\la^{j+1} \int_{t_j}^{r_{j}} u_j(t)\, dt = \frac{s-1-\eps}2
j\log\la.
\end{equation}
\end{lemma}

\begin{proof}
Note that $r_j$ is again unique. Let us assume that the statement is
not true and consider large enough $j_1$ so that $u_{j-2}(r_j)>0$
for all $j>j_1$. This also means that
\begin{equation} \lb{5.16}
u_{j-1}(t), u_j(t)>0 \text{ for } t\in I_j
\end{equation}
because $r_{j+1},r_{j+2}> t_j$.

We have $u_{j+1}(t)\le \la^{-j-1}$ for $t\in I_j$ and so by
\eqref{O} and \eqref{5.16}
\begin{equation} \lb{5.17}
u_j(t)\ge u_j(t_j) - \la^{j+1}\int_{I_j} \la^{-2j-2}\, dt  \ge
\la^{-j-1} ( \la- |I_j| )\ge \la^{-j-1}
\end{equation}
for $t\in I_j$ when $j$ is large. This, \eqref{5.4}, \eqref{5.7}
and \eqref{5.9} give
\begin{align}
\la^{j} \int_{I_j} u_{j-1}(t)\, dt & = \log
\frac{u_j(t_{j+1})}{u_j(t_j)} + \la^{j+1}
\int_{I_j}\frac{u_{j+1}(t)^2}{u_j(t)}\, dt  \notag
\\ & \le j\log\la + \la^{j+1} |I_j| \la^{-j-1}
\lb{5.17a}
\end{align}
and so
\begin{equation} \lb{5.18}
\la^{j+1} \int_{I_j} u_{j-1}(t)\, dt \le (\la+\eps)j\log\la
\end{equation}
if $j$ is large. We conclude from \eqref{5.15} and \eqref{5.18} that
there exists $a_j<r_j$, the first time in $I_j$ such that
\begin{equation} \lb{5.19}
\frac {u_j(a_j)}{u_{j-1}(a_j)} \ge \frac {s-1-\eps}{4(\la+\eps)}.
\end{equation}
Of course, sharp inequality can possibly hold only if $a_j=t_j$.
Moreover, this choice of $a_j$ and \eqref{5.18} ensure that
\[
\la^{j+1} \int_{t_j}^{a_{j}} u_j(t)\, dt \le \frac{s-1-\eps}4
j\log\la
\]
and hence by \eqref{5.15},
\begin{equation} \lb{5.20}
\la^{j+1} \int_{a_j}^{r_{j}} u_j(t)\, dt \ge \frac{s-1-\eps}4
j\log\la.
\end{equation}

Now $u_{j-2}(r_j)>0$ and \eqref{5.6} show that $E_{j-1}$ is
increasing on $[t_j,r_j]$, so that
\[
\la^{-2j}= u_j(t_j)^2\le E_{j-1}(t_j)\le E_{j-1}(a_j).
\]
Notice also that for $t\le r_j$
\[
\log \frac{u_{j+1}(t)}{u_{j+1}(t_j)} \le \la^{j+1} \int_{t_j}^t
u_{j}(\tau)\, d\tau \le \frac{s-1-\eps}2 j\log\la,
\]
which together with \eqref{5.12} gives for $t\le r_j$
\[
u_{j+1}(t)\le \la^{-\tfrac {s+1-\eps}2j}.
\]
Therefore
\begin{equation} \lb{5.20a}
E_{j+1}(a_j)\le \la^{-(s+1-\eps)j} + \sum_{l=j+2}^\infty
\la^{-2(s-\eps)l} \le \la^{-(s+1-2\eps)j} \le \la^{-2j-1}
\end{equation}
if $j$ is large. From this we have $E_{j+1}(a_j)\le \la^{-1}
E_{j-1}(t_j)$, and we obtain
\[
u_{j-1}(a_j)^2 + u_j(a_j)^2 = E_{j-1}(a_j)-E_{j+1}(a_j) \ge \tfrac
{\la -1}\la E_{j-1}(t_j) \ge \tfrac {\la -1}\la u_{j-1}(t_j)^2.
\]
This and \eqref{5.19} imply that with $c_1\equiv \tfrac{\la-1}\la
[(\tfrac{4(\la+\eps)}{s-1-\eps})^2 +1]^{-1}$,
\begin{equation} \lb{5.20b}
u_j(a_j)^2 \ge c_1 u_{j-1}(t_j)^2.
\end{equation}
Similarly as in \eqref{5.17a}, this in turn gives for $C_1\equiv
-\tfrac 12 \log c_1 + 1$
\begin{equation} \lb{5.21}
\la^{j} \int_{a_j}^{r_j} u_{j-1}(t)\, dt \le \log
\frac{u_j(r_{j})}{u_j(a_j)} + |I_j| \le \log
\frac{u_j(r_{j})}{u_{j-1}(t_j)} + C_1.
\end{equation}

Next, we claim that for large enough $j$
\begin{equation} \lb{5.21a}
u_j(r_j)\ge  (1- |I_j|) ( 1- (C_2j)^{-2}) u_{j-1}(t_j)
\end{equation}
with $C_2$ defined in \eqref{5.22} below. Assume this is not true.
Note that then $u_j(r_j)\le u_{j-1}(t_j)$, and so \eqref{5.20} and
\eqref{5.21} show that there is $b_j\in[a_j, r_j]$ such that
\begin{equation} \lb{5.22}
\frac {u_j(b_j)}{u_{j-1}(b_j)} \ge \frac {s-1-\eps}{4\la C_1} \,j
\log\la \equiv C_2j.
\end{equation}
This improves \eqref{5.19} by a factor of $j$. We now run the same
energy argument as above, with $a_j$ replaced by $b_j$ (and ignoring
the last inequality in \eqref{5.20a}), to obtain $E_{j+1}(b_j)\le
\la^{-(s-1-2\eps)j} E_{j-1}(t_j)$ and
\[
u_{j-1}(b_j)^2 + u_j(b_j)^2 \ge (1-\la^{-(s-1-2\eps)j})
u_{j-1}(t_j)^2 \ge (1-(C_2j)^{-2})u_{j-1}(t_j)^2
\]
for large $j$. \eqref{5.22} now gives $u_j(b_j)\ge
(1-(C_2j)^{-2})u_{j-1}(t_j)$, using that $(1+(C_2j)^{-2})^{-1}\ge
(1-(C_2j)^{-2})$. But then, as in \eqref{5.17}, we obtain
\begin{equation} \lb{5.22a}
u_j(r_j)\ge u_j(b_j) - \la^{j+1}\int_{b_j}^{r_j} \la^{-2j-2}\, dt
\ge ( 1- |I_j|) u_j(b_j) ,
\end{equation}
where the last inequality follows from \eqref{5.17} with $t=b_j$.
This shows \eqref{5.21a} for large enough $j$. Using \eqref{5.22a}
again, with $r_j, b_j$ replaced by $t_{j+1}, r_j$, we obtain
\begin{equation*} 
u_j(t_{j+1})\ge  (1- |I_j|)^2 ( 1- (C_2j)^{-2})u_{j-1}(t_j).
\end{equation*}

Since $\prod_{j>j_1}(1- |I_j|)^2 ( 1- (C_2j)^{-2})>0$ for large
enough $j_1$, this means that there is $c_2>0$ such that for all
large enough $j$ we have $u_j(t_{j+1})\ge c_2$. Moreover,
\eqref{5.20b} and an argument as in \eqref{5.22a} show that we
actually have for large $j$ and any $t\in[a_j,t_{j+1}]$,
\begin{equation} \lb{5.24}
c_2\le u_j(t)\le 1
\end{equation}
(with a new $c_2>0$). Then \eqref{5.21} gives
\begin{equation*} 
\la^{j+1} \int_{a_j}^{r_j} u_{j-1}(t)\, dt \le C_3
\end{equation*}
for $C_3\equiv C_1-\log c_2$. This and \eqref{5.20} means that there
is $d_j\in[a_j,r_j]$ such that
\begin{equation} \lb{5.26}
\frac {u_j(d_j)}{u_{j-1}(d_j)} \ge \frac {s-1-\eps}{8 C_3} \,j
\log\la \equiv c_3 j
\end{equation}
and
\[
\la^{j+1} \int_{a_j}^{d_{j}} u_j(t)\, dt \le \frac{s-1-\eps}8
j\log\la,
\]
and so
\[
\la^{j+1} \int_{d_j}^{r_{j}} u_j(t)\, dt \ge \frac{s-1-\eps}8
j\log\la.
\]
Thus $r_j-d_j\ge c_4j\la^{-j}$ for $c_4\equiv \tfrac{s-1-\eps}{8\la}
\log\la$.

Finally, by \eqref{5.24} we have on $[d_j,r_j]$,
\[
u_{j-1}'\le \la^{j-1} u_{j-1} - \la^j c_2^2
\]
with $u_{j-1}(d_j)\le u_j(d_j)(c_3j)^{-1}\le (c_3j)^{-1}$ by
\eqref{5.4} and \eqref{5.26}. But then for large enough $j$ we have
$u_{j-1}(d_j)< \tfrac 12 c_2^2$ and hence $u'_{j-1}(d_j)<0$. This
means that $u'_{j-1}<0$ and $u_{j-1}<\tfrac 12 c_2^2$ on
$[d_j,r_j]$. Therefore $u_{j-1}'\le - \tfrac 12\la^j c_2^2$ on
$[d_j,r_j]$, which implies
\[
u_{j-1}(r_j)\le u_{j-1}(d_j) - \tfrac 12\la^j c_2^2 c_4j\la^{-j} \le
\tfrac 12 c_2^2(1-c_4j)
\]
which is negative for large enough $j$. This contradicts
\eqref{5.16} and the proof is finished.
\end{proof}

\begin{proof}[Proof of Lemma \ref{L.5.1}]  
Let $j$ be as in the statement of Lemma \ref{L.5.2}. We will show
that if $j$ is large enough, then it also satisfies the statement of
Lemma \ref{L.5.1}.

We have (recall \eqref{5.14} and \eqref{5.15})
\begin{equation} \lb{5.27}
\la^{j+1} \int_{r_j}^{p_{j}} u_j(t)\, dt = \frac{s-1-\eps}4
j\log\la.
\end{equation}
We proceed by contradiction, so assume that  $u_{j-1}(p_j)>0$.
Notice that then $u_{j-1},u_j>0$ on $[t_j,p_j]$ (the latter by
\eqref{5.11a}). Hence \eqref{5.17} holds on this interval, and as in
\eqref{5.17a} and \eqref{5.18},
\begin{equation} \lb{5.28}
\la^{j+1} \int_{r_j}^{p_j} u_{j-1}(t)\, dt  \le (\la+\eps)j\log\la.
\end{equation}
Again, \eqref{5.27} and \eqref{5.28} show that there must be
$e_j\in[r_j,p_j]$ such that
\begin{equation} \lb{5.29}
\frac {u_j(e_j)}{u_{j-1}(e_j)} \ge \frac {s-1-\eps}{8(\la+\eps)}.
\end{equation}
and
\begin{equation} \lb{5.30}
\la^{j+1} \int_{e_j}^{p_{j}} u_j(t)\, dt \ge \frac{s-1-\eps}8
j\log\la.
\end{equation}
Now $u_{j-2}(r_j)\le 0$ gives $u'_{j-1}\le -\la^ju_j^2$ on
$[r_j,p_j]$, that is,
\begin{align}
u_{j-1}(e_j) - u_{j-1}(p_j) & \ge \la^j \int_{e_j}^{p_j} u_j(t)^2 dt
\notag
\\ & \ge \frac {\la^j} {p_j-e_j}\bigg( \int_{e_j}^{p_j} u_j(t) dt
\bigg)^2 \notag
\\ &\ge \frac{s-1-\eps}{8\la(p_j-e_j)} j\log\la \int_{e_j}^{p_j} u_j(t)
dt \lb{5.30aa}
\end{align}
by \eqref{5.30}. Since $u_{j-1}, u_j>0$ on $[e_j,p_j]$, a similar
computation as in \eqref{5.22a} shows that for $t$ in this interval
$u_j(t)\ge \tfrac 12 u_j(e_j)$ if $j$ is large. But that,
\eqref{5.29}, and \eqref{5.30aa} yield
\[
u_{j-1}(e_j) - u_{j-1}(p_j)\ge c_5 j u_{j-1}(e_j)
\]
for $c_5= \tfrac{(s-1-\eps)^2}{128\la(\la+\eps)}\log\la$ and all
large $j$. Once $j>c_5^{-1}$, this contradicts $u_{j-1}(p_j)> 0$.

Thus we have showed that if $u_{j-2}(r_j)\le 0$ and $j$ is large
enough, then $u_{j-1}(p_j)\le 0$. But then also $u_{j-1}(r_{j+1})\le
0$ because $p_j<t_{j+1}<r_{j+1}$. Lemma~\ref{L.5.1} and induction
finish the proof.
\end{proof}

Hence we have narrowed the possibility of a blowup to a
scenario where for all large $j$ there is (a single) $q_j\in
[t_j,p_j]$ such that
\begin{equation} \lb{5.30a}
u_{j-1}(q_j)=0.
\end{equation}
 That is, $u_{j-1}$ vanishes while $u_{j+1}$ is still
relatively small. Indeed, \eqref{5.14} shows that
\[
\log \frac{u_{j+1}(t)}{u_{j+1}(t_j)} \le \frac{3(s-1-\eps)}4
j\log\la,
\]
for $t\le q_j$ which together with \eqref{5.12} gives for $t\le q_j$
\begin{equation} \lb{5.31}
u_{j+1}(t) \le \la^{-(1+\tfrac{s-1-\eps}4) j} \le \la^{-(1+\eps) j}.
\end{equation}

Of course, blowup can now come only from large $u_j(q_j)$ because
all the other modes are controlled by $\la^{-j}$. Yet since
$u_{j-1}$ becomes negative on $I_{j},$ we can expect that a portion
of $u_j$ energy will be passed to the lower modes, rather than
transferred to $u_{j+1},$ making blowup unlikely. This intuition
will be confirmed if we prove the following lemma.


\begin{lemma} \lb{L.5.3}
For all large enough $j$ we have $u_{j}(q_j)\le \la^{-2}
u_{j-2}(q_{j-2})$ or $u_{j+1}(q_{j+1})\le \la^{-3}
u_{j-2}(q_{j-2})$.
\end{lemma}

Let us first complete the proof of Theorem \ref{thm3} given this
lemma.

\begin{proof}[Proof of Theorem \ref{thm3}] Choose large enough $j_1$ and set $C\equiv
\la^{j_1}u_{j_1}(q_{j_1})$. Then the lemma and induction show that
there is a sequence $j_l\to\infty$ such that $\la^{j_l}
u_{j_l}(q_{j_l})\le C$. Since on $[q_{j_l},t_{j_l+1}]$ we have
$u_{j_l}'\le 0$, we also have there $\la^{j_l} u_{j_l}\le C$ which
gives $u_{j_l+1}'\le C\la u_{j_l+1}$. But then \eqref{5.7} and
\eqref{5.31} show
\[
(\eps j_l-1)\log\la\le \log \frac
{u_{j_l+1}(t_{j_l+1})}{u_{j_l+1}(q_{j_l})} \le
C\la(t_{j_l+1}-q_{j_l})\le C\la.
\]
This is a contradiction when $l$ is large.
\end{proof}

Thus, we are left with proving the lemma.

\begin{proof}[Proof of Lemma \ref{L.5.3}]
Notice that \eqref{5.17} holds for $t\le q_j$ and so
\begin{equation} \lb{5.32}
u_{j}(q_j) \ge \la^{-j-1}.
\end{equation}
Also $u_j$ obviously decreases on $[q_j,q_{j+1}]$. Let now $j$ be
large enough so that Lemma \ref{L.5.1} holds for any $j'\ge j-2$ in
place of $j$. In particular, $u_{j'-1}(q_{j'})=0$ with $q_{j'}$
defined above. Let us denote
\begin{equation} \lb{5.32a}
B\equiv u_{j-2}(q_{j-2}) \ge \la^{-j+1}.
\end{equation}
By \eqref{5.12}, \eqref{5.31}, and \eqref{5.32a}
we have $E_{j-2}(q_{j-2})\le \tfrac {25}{16} B^2$ if $j$ is large
enough. Since $u_{j-3}\le 0$ on $[q_{j-2},T)$, \eqref{5.6} gives
$E_{l}(t)\le \tfrac {25}{16} B^2$ for $l\ge j-2$ and $t\ge
q_{j-2}$, in particular,
\begin{equation} \lb{5.34}
u_{l}(t)\le \tfrac 54\, B \text{ for } l\ge j-2 \text{ and } t\ge
q_{j-2}.
\end{equation}

Let us again proceed by contradiction and assume that
\begin{equation} \lb{5.35}
u_{j}(q_j)> \tfrac B {\la^{2}} \text{ and } u_{j+1}(q_{j+1})>
\tfrac B {\la^{3}}.
\end{equation}
We define $f_j\in[q_j, q_{j+1}]$ to be the first time such that
\begin{equation} \lb{5.36}
u_{j}(f_j)=u_{j+1}(f_j)
\end{equation}
(recall that $u_j(q_j)\ge\la^{-j-1}\ge u_{j+1}(q_j)$ and
$u_j(q_{j+1})=0\le u_{j+1}(q_{j+1})$). Then we must have for any
$t\in[q_j,f_j]$,
\begin{equation} \lb{5.37}
u_{j}(t)\ge \tfrac B{2\la^3} \quad \big(\ge \tfrac 12 \la^{-j-2}
\big)
\end{equation}
because otherwise \eqref{5.6}, \eqref{5.12}, \eqref{5.31}, the
definition of $f_j$, and $u_{j-1}(t)\le 0$ show that
\[
E_{j+1}(q_{j+1})= E_j(q_{j+1})\le E_j(t)\le 3u_j(t)^2\le \tfrac
{B^2}{\la^{6}},
\]
which would mean $u_{j+1}(q_{j+1})\le B\la^{-3}$, contradicting the
assumption. We assume here again that $j$ is large enough, so that
$E_{j+2}(t) \leq u_j(t)^2$ for $t \in [q_j, f_j].$ Therefore
\eqref{5.37} holds and there is a first time $g_j\in[q_j, f_j]$ such
that $u_{j+1}(g_j) = 2 \la^{-(1+\eps) j}$. From \eqref{O},
\eqref{5.31}, and $u_j(t)\ge u_j(f_j)\ge B(2\la^3)^{-1}$ and
$u_{j+1}(t)\le u_{j+1}(g_j)=2 \la^{-(1+\eps) j}$ for $t\in[q_j,g_j]$
we get
\begin{align}
\la^{j+1}\int_{q_j}^{g_j} u_j(t)^2 \,dt  & \ge \frac B{4
\la^{3-(1+\eps) j}} \la^{j+1}\int_{q_j}^{g_j} u_j(t) u_{j+1}(t)
\,dt \notag
\\ & \ge \frac B{4 \la^{3-(1+\eps)
j}} (u_{j+1}(g_j) - u_{j+1}(q_j)) \notag
\\ & \ge \frac B{4\la^3}. \lb{5.38}
\end{align}

Next, we will show that there is $h_j\in[q_j,g_j]$ such that
\begin{equation} \lb{5.39}
u_{j-1}(h_j)\le -\tfrac B{10\la^5}, \quad u_{j}(h_j)\ge \tfrac
B{2\la^3}, \quad u_{j+1}(h_j)\le 2\la^{-(1+\eps)j}.
\end{equation}
The second and third inequality are automatic when $h_j\le g_j$
(by \eqref{5.37} and the definition of $g_j$), so let us assume
that for all $t\in[q_j,g_j]$ we have $- B(10\la^5)^{-1} <
u_{j-1}(t) (\le 0)$. Then on $[q_j,g_j]$ we have by $u_j\ge
B(2\la^3)^{-1}$, \eqref{O} and \eqref{5.34},
\[
u_{j-1}'\le \la^{j-1} \tfrac B{10\la^5}\, \tfrac {5B}4  -\la^j
u_j^2\le -\tfrac 12 \la^j u_j^2.
\]
Hence from \eqref{5.38},
\[
u_{j-1}(g_j)\le u_{j-1}(q_j) - \tfrac B{8\la^4} = - \tfrac
B{8\la^4}\le - \tfrac B{10\la^5},
\]
a contradiction with the assumption.

Therefore \eqref{5.39} holds for some $h_j\in[q_j,g_j]$. Moreover,
$u_{j-1}\le -B(10\la^5)^{-1}$ on $[h_j,f_j]$ because whenever
equality holds, then by \eqref{5.34} and \eqref{5.37},
\[
u_{j-1}'\le \la^{j-1}\,  \tfrac B{10\la^5} \tfrac {5B}4 - \la^j
\Big(\tfrac B{2\la^3} \Big)^2 < 0.
\]
Therefore by \eqref{5.34}, on $[h_j,f_j]$
\[
\frac {u_j'}{u_j} \le \la^j u_{j-1} \le -\la^j\, \frac B{10\la^5}
\le -\la^j \, \frac {\tfrac 45 u_j}{10\la^5} = -\frac
2{25\la^6}\la^{j+1} u_j.
\]
From this,
\[
\log \frac{u_j(f_j)}{u_j(h_j)} \le -\frac 2{25\la^6}
\la^{j+1}\int_{h_j}^{f_j} u_j(t)\, dt ,
\]
which together with \eqref{5.34} and \eqref{5.37} shows that
\begin{equation*} 
\la^{j+1}\int_{h_j}^{f_j} u_j(t)\, dt \le \frac {25\la^6}2 \log
\frac{5\la^3}2 \equiv C_6.
\end{equation*}
But then \eqref{5.32a}, \eqref{5.36}, \eqref{5.37}, \eqref{5.39},
and \eqref{O} yield
\[
\log \frac{\la^{\eps j}}{4\la^2}\le  \log \frac
B{4\la^3\la^{-(1+\eps)j}} \le \log \frac
{u_{j+1}(f_j)}{u_{j+1}(h_j)} \le \la^{j+1}\int_{h_j}^{f_j}
u_j(t)\, dt \le C_6,
\]
a contradiction when $j$ is large. The lemma is proved.
\end{proof}

\section{Appendix}

Here, for the sake of completeness, we sketch the argument giving
the local existence of solutions in the Obukhov and KP models. These
two cases (as well as any combination) are handled identically; for
simplicity we will consider the Obukhov model. We will also look at
a more general branching case, since the result extends naturally
and without extra effort. Thus, we look at a largest dyadic cube
$Q^0$ of generation zero, and assume it has $d$ children $Q^1_l$
belonging to the first generation. Each cube $Q$ of generation $j$
has in its turn $d$ children of generation $j+1.$ Given a cube $Q$
of generation $j,$ we denote $\tilde{Q}$ its unique parent and
$C^{1}(Q)$ the set of its $d$ children. Likewise, we denote $C^k(Q)$
the set of all descendants of $Q$ of generation $j+k.$ Let us denote
$j(Q)$ the generation of any given cube $Q$ in our branching tree.
The branched Obukhov model is given by the following system of
differential equations
\begin{equation}\label{Obbra}
\frac{d}{dt} u_Q = \lambda^{j(Q)} u_{\tilde{Q}} u_Q -
\lambda^{j(Q)+1} \sum\limits_{Q' \in C^1(Q)} u_{Q'}^2
\end{equation}
for each $Q$, with $u_{\tilde{Q}}\equiv 0$ when $Q=Q^0$. We say that
$U \equiv \{u_Q\}$ belongs to the Sobolev space $H^s$ if
\[ \|U\|_{H^s}^2 \equiv \sum\limits_{Q} \lambda^{2sj(Q)}
|u_Q|^2 < \infty. \]

Consider an equivalent integral equation reformulation of
\eqref{Obbra}, given by
\begin{equation}\label{Obbraint}
u_Q(t) = u_Q(0)+ \int_0^t \left(\lambda^j u_{\tilde{Q}}(\tau)
u_Q(\tau) - \lambda^{j+1} \sum\limits_{Q' \in C^1(Q)} u_{Q'}(\tau)^2
\right)\,d\tau.
\end{equation}
Recall one version of the well-known Picard's fixed point theorem.
\begin{theorem}[Picard]\label{picard}
Let $X$ be a Banach space and $\Gamma$ a bilinear operator $\Gamma:
\,X\times X \mapsto X$ such that for any $U,V \in X$ we have
\begin{equation}\label{gbound}
\|\Gamma(U,V)\|_X \leq \eta \|U\|_X \|V\|_X.
\end{equation}
Then for any $U_0 \in X$ satisfying $4\eta \|U_0\|_X <1$ the
equation $U = U_0 + \Gamma(U,U)$ has a unique solution $U \in X$
such that $\|U\|_X \leq 1/2\eta.$
\end{theorem}
Using this theorem we are going to prove
\begin{theorem}\label{locex}
Given any $\{u_Q(0)\} \in H^s,$ $s \geq 1,$ there exists
$T=T(\|u_Q\|_{H^s})>0$ such that there is a unique solution $u_Q(t)$
of the branching Obukhov system \eqref{Obbraint} which belongs to
$C([0,T], H^s).$
\end{theorem}
\begin{proof}
Let us define $U_0(t)\equiv \{u_Q(0)\}$ for all $t$ and
\[ \gamma(U,V)_Q(t) = \lambda^{j(Q)} u_{\tilde{Q}}(t) v_Q(t)
- \lambda^{j(Q)+1} \sum\limits_{Q' \in C^1(Q)} u_{Q'}(t) v_{Q'}(t) ,
\] and
\[ \Gamma(U,V)_Q(T) = \int\limits_0^T \gamma(U,V)_Q(t) \,dt.\]
The result will follow from Picard's theorem if we verify the bound
\eqref{gbound} for $\Gamma.$ If $s \geq 1,$ we have
\begin{align*}
\|\gamma(U,V)(t)\|_{H^s}^2 & = \sum\limits_Q \lambda^{2sj(Q)}
\left(\lambda^{j(Q)} u_{\tilde{Q}}(t) v_Q(t) - \lambda^{j(Q)+1}
\sum\limits_{Q' \in C^1(Q)} u_{Q'}(t) v_{Q'}(t)\right)^2
\\ & \leq \sum\limits_Q \lambda^{2sj(Q)} (d+1)
\left(\lambda^{2j(Q)} u_{\tilde{Q}}(t)^2 v_Q(t)^2 +
\lambda^{2j(Q)+2} \sum\limits_{Q' \in C^1(Q)} u_{Q'}(t)^2
v_{Q'}(t)^2\right)
\\ & \leq (d+1) \lambda^{2s}\|U(t)\|_{H^s}^2 \left( \sum\limits_Q \lambda^{2j(Q)}
v_Q(t)^2 +  \sum\limits_Q \sum\limits_{Q' \in C^1(Q)}
\lambda^{2j(Q)+2} v_{Q'}(t)^2 \right)
\\ & \leq
2(d+1)\lambda^{2s}\|U(t)\|_{H^s}^2\|V(t)\|_{H^s}^2.
\end{align*}
Then
\[ \|\Gamma(U,V)\|_{C([0,T], H^s)} \leq C(d,\lambda)
\int\limits_0^T \|U(t)\|_{H^s} \|V(t)\|_{H^s} \,dt \leq
C(d,\lambda)T \|U\|_{C([0,T], H^s)}\|V\|_{C([0,T], H^s)}. \]
Choosing a small enough $T>0$ completes the proof.
\end{proof}

\smallskip

\noindent {\bf Acknowledgement.}\rm The work of AK has been
supported in part by the Alfred P. Sloan Research Fellowship and
NSF-DMS grant 0314129. The work of AZ has been supported in part by
the NSF-DMS grant 0314129.

\end{document}